\documentclass[12pt]{article}
\usepackage{amssymb,amsmath,amscd}

\newtheorem{proposition}{Proposition}
\newtheorem{lemma}[proposition]{Lemma}
\newtheorem{remark}[proposition]{Remark}

\def\A{{\cal A}}
\def\Ad{{\rm Ad}}
\def\doubleoplus#1#2{{\bigoplus_{\stackrel{#1}{\scriptscriptstyle#2}}}\ }
\def\Co{\mathbb{C}} 
\def\cotens#1#2#3{{#1\,\square_{#2}\,#3}}

\def\fidi{\hskip5pt \vrule height4pt width4pt depth0pt \par}

\def\K{{\cal K}}

\def\H{\mathbb{H}} 
\def\id{{\rm id}}
\def\L{{\cal L}}
\def\LieG{\mathfrak{g}}

\def\Na{\mathbb{ N}}  
\def\Pol{{\rm Pol}}
\def\R{{\cal R}}
\def\Re{\mathbb{R}}  

\def\sphere{\mathbb{S}}
\def\sph4{\Sigma_q^4}
\def\Tr{{\rm Tr}}
\def\tr{{\rm tr}}
\def\u4{{\bf u}(4)}

\begin{document}
\title{ Noncommutative Instantons on the 4--Sphere \\from Quantum Groups}

\author{F.Bonechi${}^{1,2}$, N.Ciccoli${}^3$, M.Tarlini${}^{1,2}$}

\date{December 2000}

\maketitle

\centerline{{\small  ${ }^1$ INFN Sezione di Firenze}}
\centerline{{\small ${ }^2$ Dipartimento di Fisica, Universit\`a
di Firenze, Italy. }} \centerline{{\small ${ }^3$ Dipartimento di
Matematica, Universit\`a di Perugia, Italy. }} \centerline{{\small
e-mail: bonechi@fi.infn.it, ciccoli@dipmat.unipg.it,
tarlini@fi.infn.it}}

\begin{abstract}
\noindent
{We describe an approach to the noncommutative instantons on the 4--sphere
based on quantum group theory. We quantize the Hopf bundle $\sphere^7\to
\sphere^4$ making use of the concept of quantum coisotropic subgroups.
The analysis of the semiclassical Poisson--Lie structure of $U(4)$ shows
that the diagonal $SU(2)$ must be conjugated to be properly quantized.
The quantum coisotropic subgroup we obtain is the standard $SU_q(2)$; it
determines a new deformation of the 4--sphere $\sph4$ as the
algebra of coinvariants in $\sphere_q^7$.
We show that the quantum vector bundle associated to the fundamental
corepresentation of $SU_q(2)$ is finitely generated and projective and
we compute the explicit projector. We give the unitary representations
of $\sph4$, we define two 0--summable Fredholm modules and we compute the
Chern--Connes pairing between the projector and their characters.
It comes out that even the zero class in cyclic homology is non trivial.}
\end{abstract}

\thispagestyle{empty}

\section{Introduction}
Since the work \cite{NS} on instantons on noncommutative $\Re^4$ a
lot of attention has been devoted to the problem of gauge theories
on noncommutative four manifolds.
In ordinary differential geometry, the topological properties of instantons
in $\Re^4$ are better understood by studying fibre bundles on the sphere
$\sphere^4$.
In noncommutative geometry this is not an easy task: it is more
natural to define the problem directly on the noncommutative
sphere.

Very recently, in \cite{ConnesLandi} and \cite{DLM} two different deformations
of $\sphere^4$ were proposed.
The one in \cite{ConnesLandi} preserves the property of having zero the first
Chern class which is not trivial in \cite{DLM}.
In this second case the deformation is a suspension of the quantum
3--sphere $SU_q(2)$ obtained by adding a central generator.

In this paper we propose an alternative approach, based more directly on
quantum groups and on Hopf algebraic techniques.

In noncommutative geometry finitely generated projective modules, {\it i.e.}
the quantum vector bundles, are the central object to develop gauge theories.
From this point of view  there is no obvious notion of structure group.
Quantum groups provide a construction of quantum vector bundles
which is closer to ordinary differential geometry.
The first attempts go back to \cite{Durd}, \cite{Pflaum} and \cite{BrMj},
where the gauge theory is developed starting from the notion of Hopf--Galois
extension, which is the analogue of principal bundles in the Hopf algebra
setting, see \cite{Schneider}.
The associated quantum vector bundles have a Hopf algebra
on the fiber and, if they admit a connection, are finitely generated and
projective modules \cite{DabGroHaj}.

Although this definition works in principle, it is not enough to explain
all known interesting examples.
This problem is better understood if we concentrate on the specific class
of principal bundles given by homogeneous spaces.
A quantum homogeneous space is an example of Hopf--Galois extension
only if it is obtained as quotient by a quantum subgroup ({\it i.e.} a Hopf
algebra quotient).
But quantum subgroups are very rare. For instance between the  quantum
$2$--spheres introduced by Podle\`s in \cite{Podles}
only one, the standard one, is such an example.
It is necessary to generalize the notion of subgroup, allowing a more general
quotient procedure.
This is possible by using {\it quantum coisotropic subgroups}: they are
quotient by a coideal, right (or left) ideal, so that they inherit only
the coalgebra, while the algebra structure is weakened to a right (or left)
module.
Their semiclassical interpretation is illuminating: in a Poisson--Lie group
every Poisson (resp. coisotropic) subgroup can be quantized to a quantum
(resp. coisotropic) subgroup (see \cite{Ciccoli}).
Nevertheless conjugation, which does not change topology, can break Poisson
properties: for instance a subgroup conjugated to a Poisson subgroup can be
only coisotropic or can have no Poisson properties at all (see for instance
$SL(2,\Re)$ in \cite{BCGST2}).
Coisotropic subgroups can be quantized and give rise to inequivalent quantum
homogeneous spaces: for instance all the Podle\`s quantum spheres are obtained
as quotient of coisotropic $U(1)$. The general scheme to describe such examples
could be the so called $C$--Galois extensions, see for instance
\cite{BrzHaj,Brz}.

The principal bundle on $\sphere^4$ corresponding to $SU(2)$ instantons
with charge $-1$ has $\sphere^7=U(3)\backslash U(4)$ as total space and
the action on the fibre is obtained by considering $SU(2)$ as diagonal
subgroup of $U(4)$.
In this description $\sphere^4$ is the double coset $U(3)\backslash
U(4)/SU(2)$.
In the quantum setting, odd spheres were obtained in \cite{VS}
as homogeneous spaces of $U_q(N)$ with respect to the quantum subgroup
$U_q(N-1)$ so that the left quotient is easily quantized.
The right quotient is more problematic because the diagonal
$SU(2)$ doesn't survive in the quantization of $U(4)$; indeed the analysis
of the limit Poisson structure on $U(4)$ shows that it is not coisotropic.
We then have to look for coisotropic subgroups in
the conjugacy class of the diagonal one.
It comes out that there is at least one which defines what we call the
Poisson Hopf bundle in $\sphere^4$ (Proposition \ref{su2dg}).
In this bundle, which is topologically equivalent to the usual Hopf bundle,
both the total and the base spaces are Poisson manifolds and the projection
is a Poisson map.
Its quantization is straightforward: the quantum coisotropic subgroup turns
out to be equivalent as coalgebra to $SU_q(2)$ (Proposition \ref{suq2})
and the algebra of functions over the quantum 4--sphere $\sph4$ is then
obtained as the subalgebra of coinvariants in $\sphere^7_q$ with respect
to this $SU_q(2)$ (Proposition \ref{s4q}).
This deformation of the algebra of functions on $\sphere^4$ is different
from those introduced in \cite{ConnesLandi} and \cite{DLM}.
We then study the quantum vector bundle associated to the
fundamental corepresentation of $SU_q(2)$ and give the explicit projector
(Proposition \ref{projector}).
We describe the unitary irreducible representations of $\sph4$
(Equation \ref{rep0} and \ref{irrep}); there is a 1--dimensional
representation and an infinite dimensional one realized by trace class
operators (Proposition \ref{trace_class}).
Finally we study the Chern class in cyclic homology of the projector and
compute the Chern--Connes pairing with a trace induced by the trace class
representation (Proposition \ref{chern_class}). It comes out that, on the
contrary with \cite{ConnesLandi} and \cite{DLM}, they are all non trivial.
This result is the analogue of what was obtained in
\cite{MasNakWat,HajacMajid}  for the standard  Podle\`s 2--sphere.


\section{Quantization of Coisotropic Subgroups} A Poisson--Lie
group $({\cal G},\{,\})$ is a Lie group ${\cal G}$ with a Poisson
bracket $\{,\}$ such that the multiplication map $m:{\cal G}\times
{\cal G}\to {\cal G}$ is a Poisson map with respect to the
product Poisson structure in ${\cal G}\times {\cal G}$.
The Poisson bracket $\{,\}$ is identified by a bivector $\omega$
({\it i.e.} a section of $\bigwedge^2T{\cal G}$) such that
$\{\phi_1,\phi_2\} (x)=\omega(x)(d_x\phi_1,d_x\phi_2)$.
(For more details see \cite{Chari_Pressley} and \cite{Vai})

Every Poisson--Lie group induces a natural bialgebra structure on
$\LieG= {\rm Lie}({\cal G})$ which will be called the tangent bialgebra
of ${\cal G}$.
Indeed, $\delta:\LieG\to\LieG\wedge\LieG$ is defined by
$\langle X,d_e\{f,g\}\rangle = \langle
\delta(X),f\otimes g\rangle$, where $X\in\LieG$ and $f,g\in
C^\infty({\cal G})$.

The point we want to discuss is the behaviour of subgroups and
corresponding homogeneous spaces with respect to the Poisson
structure.
A Lie subgroup ${\cal H}$ of ${\cal G}$ is called a Poisson--Lie
subgroup if it is also a Poisson submanifold of ${\cal G}$, {\it i.e.}
if the immersion map $\imath:{\cal H}\to {\cal G}$ is a Poisson map.
There are various characterizations for such subgroups: as invariant
subspaces for the dressing action or as union of symplectic leaves \cite{LuW}.

The property of being a Poisson--Lie subgroup is, evidently, a very
strong one.
We need then to characterize a family of subgroups satisfying weaker
hypothesis with respect to the Poisson structure.

In Poisson geometry a submanifold $N$ of a Poisson manifold
$(M,\omega)$ is said to be coisotropic if $\omega\big|
_{{\rm Ann}(TN)}=0$, where ${\rm Ann} (T_xN)=\{ \alpha\in T^*_x(M)\,
|\,\alpha(v)=0\quad \forall v\in T_xN\}$.
Coisotropy can be formulated very neatly as an algebraic property at the
function algebra level (see \cite{Vai}). Indeed a locally closed submanifold
$N$ of the Poisson manifold $(M,\omega)$ is coisotro\-pic if and only
if for every $f,g\in\Co^{\infty}(M)$
$$
\label{idea} f\big|_N=0, \, g\big|_N=0 \Rightarrow \{
f,g\}\big|_N=0 ~~~.
$$
Thus locally closed coisotropic submanifolds correspond to
manifolds whose defining ideal is not a Poisson ideal but only a
Poisson subalgebra.

A Lie subgroup ${\cal H}$ of $({\cal G},\omega)$ is said a {\it
coisotropic subgroup} if it is coisotropic as Poisson submanifold.
In the connected case there are nice characterizations, as shown
for example in \cite{Lu}; we will need the following one:
\medskip
\begin{proposition}
\label{coisotrop} A connected subgroup ${\cal H}$ of $({\cal
G},\omega)$ with $\mathfrak{h}={\rm Lie}{\cal H}$ is coisotropic
iff $~\delta(\mathfrak{h})\subset\mathfrak{g}\wedge\mathfrak{h}$
and it is Poisson--Lie iff  $~\delta(\mathfrak{h})\subset\mathfrak{h}
\wedge\mathfrak{h}$.
\end{proposition}
\medskip

Given a Poisson--Lie group ${\cal G}$ and a coisotropic subgroup
${\cal H}$ the natural projection map ${\cal G}\to {\cal H\backslash G}$
coinduces a Poisson structure on the quotient.
If ${\cal K}$ is a second subgroup of ${\cal G}$ a condition which
guarantees that even the projection on the double coset is Poisson is
given by the following:
\medskip
\begin{proposition}[\cite{Lu}]
\label{prop_Lu}
Let $(M,\omega_o)$ be a Poisson manifold with a Poisson action of a
Poisson--Lie group $({\cal G},\omega)$. Let ${\cal K}$ be a coisotropic
connected subgroup of ${\cal G}$.
If the orbit space $M/{\cal K}$ is a manifold there exists a unique
coinduced Poisson bracket such that the natural projection
$M \to M/{\cal K}$ is Poisson.
\end{proposition}
\medskip
We now recall how these concepts can be translated in a Hopf
algebra setting, (see \cite{BCGST,Ciccoli} for more details).
Given a real quantum group $(\A,*,\Delta,S,\epsilon)$ we will call
{\it real coisotropic quantum right $($left$)$ subgroup} $(\K,
\tau_\K)$ a coalgebra, right (left) $\A$--module $\K$ such that:
\begin{itemize}
\item[i)] there exists a surjective linear map $\pi:
\A\to \K$, which is a morphism of coalgebras and of
$\A$--modules (where $\A$ is considered as a module on itself
via multiplication);
\item[ii)] there exists an antilinear map $\tau_\K:
\K\to \K$ such that
$\tau_\K\circ\pi=\pi\circ\tau$, where $\tau = *\circ S$.
\end{itemize}
A $*$--Hopf algebra ${\cal S}$ is said to be a real quantum subgroup
of $\A$ if there exists a $*$--Hopf algebra epimorphism
$\pi:\A\to {\cal S}$; evidently this is a particular
coisotropic subgroup.
We remark that a coisotropic quantum subgroup is not in general
a $*$--coalgebra but it has only an involution $\tau_\K$ defined on it.

Right (left) coisotropic quantum subgroups are obviously
characterized by the kernel of the projection, which is a
$\tau$--invariant two--sided coideal, right (left) ideal in $\A$.
It is easy to verify that if the kernel is also $*$--invariant
then it is an ideal and the quotient is a real quantum subgroup.

\medskip
A $*$--algebra $B$ is said to be an embeddable quantum left
(right) $\A$--homo\-geneous space if there exists a coaction $\mu:
B\to B\otimes \A$, ($\mu: B\to \A\otimes B$) and an injective
morphism of $*$--algebras $\imath :B\to\A$ such that $\Delta\circ
\imath = (\imath\otimes \id)\circ\mu$ ($\Delta\circ \imath =
(\id\otimes \imath)\circ\mu$).

Embeddable quantum homogeneous spaces can be obtained as the space of
coinvariants with respect to the coaction of coisotropic quantum
subgroups.
For instance if ${\cal K}$ is a right (left) subgroup and
$\Delta_\pi=({\rm id}\otimes\pi)\Delta$ ($\,{}_\pi\Delta=
(\pi\otimes {\rm id})\Delta\,$), then
$$
B^\pi = \{a\in\A \,|\, \Delta_\pi a = a\otimes\pi(1)\}\quad(\,{}^\pi B=
\{a\in\A \,|\, {}_\pi\Delta a = \pi(1)\otimes a\}\,),
$$
is an homogeneous space with $\mu=\Delta$.

If $\rho:V\to {\cal K}\otimes V$ is a corepresentation of
${\cal K}$, we define the cotensor product as
$$
\cotens{{\cal A}}\rho{V} = \{F\in {\cal A}\otimes V \ |\
(\Delta_\pi \otimes \id)\, F =(\id\otimes\rho)\, F\}\;.
$$
We have that $\cotens{{\cal A}}\rho{V}$ is a left $B^\pi$--module.
Let $\rho$ be unitary and $\{e_i\}$ be an orthonormal basis of
$V$; if $F=\sum_i F_i\otimes e_i$, let's define $\langle
F,G\rangle = \sum_i F_i G_i^*$.
It is shown in \cite{BCGST} that $\langle,\rangle$ is a sesquilinear
form on $\cotens{{\cal A}}\rho{V}$ with values in $B^\pi$.

\medskip
The correspondence between coisotropic quantum subgroups and
embeddable quantum homogeneous spaces is bijective only provided some
faithful flatness conditions on the module and comodule structures are
satisfied (see \cite{MS} for more details).

The role of coisotropic subgroups can also be appreciated in the context
of formal and algebraic equivariant quantization.
While it is known that not every Poisson homogeneous space admits such
quantization, it holds true that every quotient of a Poisson--Lie group
by a coisotropic subgroup can be equivariantly quantized.
Although such quotients do not exhaust the class of quantizable Poisson
spaces they provide a large subclass in it.
Furthermore in functorial quantization they correspond to embeddable
quantum homogeneous spaces. More on the subject can be found in \cite{EK}.


\section{The Classical Instanton with $k=-1$}
In this section we review the construction of the principal bundle
corresponding to instantons with topological charge $k=-1$ (see
\cite{At}). We denote with $\H$ the quaternions generated
by $i$, $j$, $k$ with the usual relations $i^2=j^2=k^2=-1$, and
$ij=-ji=k$, $jk=-kj=i$, $ki=-ik=j$. The total space of the bundle
is defined as $E = \{(q_1,q_2)\in \H^2 \,|\, |q_1|^2+|q_2|^2=1\}$,
the base space is $P_1(\H)=\{[(q_1,q_2)]\,|\,
(q_1,q_2)\simeq(q_1\lambda,q_2\lambda),(q_1,q_2)\in \H^2, \lambda\in \H\}$
and the bundle projection is $p(q_1,q_2)=[(q_1,q_2)]$. The fibre is $SU(2)$
which acts on $\H^2$ by the diagonal right multiplication of
quaternions of unit modulus. The quaternionic polynomial functions
$B=\Pol(P_1(\H))$ on the base space are generated by
$R=q_1\bar{q_1}$ and $Q=q_1\bar{q_2}$, with the relation $|Q|^2=R(1-R)$.

The fundamental representation of $SU(2)$ can be realized again by
right multiplication of unit quaternions on $\H$. The space $\cal
E$ of sections of the associated vector bundle is the space of
equivariant functions $s:E\to \H$, {\it i.e.} such that
$s(q_1,q_2)\lambda=s(q_1\lambda, q_2\lambda)$, for $|\lambda|=1$.
It is generated as a left $B$--module by $s_1(q_1,q_2)=q_1$ and
$s_2(q_1,q_2)=q_2$ and it has an hermitian structure
$\langle,\rangle: {\cal E}\times {\cal E}\to B$ given by $\langle
s_1,s_2\rangle = s_1\bar{s_2}$.

We can define $G\in M_2(B)$ with $G_{ij}=\langle
s_i,s_j\rangle$. By direct computation we obtain that
\begin{equation}
\label{projector_Q}
G=G^2=\left(
\begin{array}{cc}
R & Q \cr \bar{Q} & 1-R
\end{array}
\right) \;.
\end{equation}
It is easy then to verify that ${\cal E}\simeq B^2\,G$.

\medskip
For our future purposes, we have to describe this bundle in a Hopf
algebraic language. We first remark that $E$ is isomorphic to
$\sphere^7=U(3)\backslash U(4)$ and $P_1(\H)$ to
$\sphere^4=U(3)\backslash U(4)/SU(2)$.

Let $t_f=\{t_{ij}\}_{ij=1}^4$ define the fundamental representation
of $U(4)$. Then $\Delta(t_{ij})\!\!= \sum_k t_{ik}\otimes t_{kj} $
and let $\ell:\Pol(U(4))\to\Pol(U(3))$ be the Hopf algebra
projection defined by $\ell(t_{4j})=\ell(t_{j4})=0$ for $j=1,2,3$,
and $\ell(t_{44})=1$. The algebra of polynomial functions on
$\sphere^7$ is given by the coinvariants ${}^{\ell}\Pol(U(4))$ and
it is generated by $z_i=t_{4i}$ with the relation $\sum_i
|z_i|^2=1$. Let $r:\Pol(U(4))\to\Pol(SU(2))$ be the Hopf algebra
projection defined by
$$
r(t) = \left(
\begin{array}{cccc}
\alpha & \beta & 0 & 0 \cr -\beta^* & \alpha^* & 0 & 0\cr 0 & 0 &
\alpha & \beta \cr 0 & 0 & -\beta^* & \alpha^*
\end{array}
\right)\,,~~~~~~~~ |\alpha|^2+|\beta|^2=1 \;.
$$
As usual $\Pol(U(4)/SU(2))$ is obtained as the space of
coinvariants $\Pol(U(4))^r$. The algebra of polynomial functions on
$\sphere^4$ is ${}^\ell\Pol(U(4))\cap \Pol(U(4))^r$ and is generated by
$R=|z_1|^2+|z_2|^2$, $A=z_1z_3^*+z_2z_4^*$ and $B=z_1z_4-z_2z_3$,
with the relation $|A|^2+|B|^2=R(1-R)$.

Let $\tau_f: \Co^2\to \Pol(SU(2))\otimes\Co^2$ be the fundamental
corepresentation of $\Pol(SU(2))$
$$
\tau_f\left(\begin{array}{c}e_1\cr e_2\end{array}\right)=
\left(\begin{array}{cc}\alpha&\beta\cr
-\beta^*&\alpha\end{array}\right)\otimes
\left(\begin{array}{c}e_1\cr e_2\end{array}\right)\;.
$$
The left $\Pol(\sphere^4)$--module of sections of the associated vector
bundle is obtained as ${\cal E}=\cotens{\Pol(\sphere^7)}{\tau_f}\Co^2$.
As a $\Pol(\sphere^4)$--module, ${\cal E}$ is generated by
$$
f_1=\left(\begin{array}{c}z_1\cr z_2\end{array}\right)\,,~~
f_2=\left(\begin{array}{c}z_2^*\cr-z_1^*\end{array}\right)\,,~~
f_3=\left(\begin{array}{c}z_3\cr z_4\end{array}\right)\,,~~
f_4=\left(\begin{array}{c}z_4^*\cr-z_3^*\end{array}\right)\,.
$$
With the usual hermitian structure we define $G\in \Pol(\sphere^4)
\otimes M_4(\Co)$ with $G_{ij}=\langle f_i,f_j\rangle$ and obtain that
\begin{equation}
\label{projector_C}
 G = G^2=\left(\begin{array}{cccc} R & 0 & A & B\cr
0& R& -B^*& A^*\cr A^* & -B & 1-R & 0\cr B^* & A & 0 & 1-R
\end{array}\right)\;.
\end{equation}

With the usual representation of $\H$ as $\Co^2$, where
$(z_1,z_2)$ is identified with $z_1+z_2j$, it is easy to verify
that $Q=A-Bj$, $f_1=q_1$, $f_2=-jq_1$, $f_3=q_2$ and $f_4=-jq_2$.
Once we introduce the representation of the quaternions with Pauli
matrices it is easy to verify that (\ref{projector_Q}) and
($\ref{projector_C}$) define the same projector.


\section{Poisson Hopf bundle on $\sphere^4$}
Let us identify $\LieG=\u4={\rm Lie}\,U(4)$ with its defining
representation by antihermitian $4\times 4$ matrices.
The $SU(2)$ generators of the Dynkin diagram are, for $i=1,2,3$
$$
H_i= i\left(e_{ii}-e_{i+1,i+1}\right) ~,~~~ E_i = \frac{1}{2i}
\left(e_{i,i+1}+ e_{i+1,i}\right) ~,~~~ F_i = \frac{1}{2}
\left(e_{i,i+1}-e_{i+1,i}\right) \;,
$$
where $e_{ij}$ are the elementary matrices with entries
$(e_{ij})_{kl}=\delta_{ik} \delta_{jl}$ and the central generator is
$H=i\mathbb{I}$.
The Poisson--Lie structure of $U(4)$ is defined by the canonical
coboundary bialgebra given on these generators by
\begin{equation}
\label{coboundary} \delta_{\Re}(H_i)=0,\quad\delta_{\Re}(H)=0,
\quad\delta_{\Re}(E_i)=E_i\wedge H_i,\quad\delta_{\Re}(F_i)=F_i\wedge H_i\;.
\end{equation}

The generators $h=\frac14 H_1+\frac12 H_2 +\frac34 H_3 +\frac34 H$
and $\{H_i,E_i,F_i\}_{i=1,2}$ define the embedding of ${\bf u}(3)$
in ${\bf u}(4)$ that we want to study; from relations (3) we have
that $\delta_{\Re}({\bf u }(3))\subset {\bf u}(3)\wedge{\bf u}(3)$
so that $U(3)$ is a Poisson Lie subgroup.

Let us fix on $\sphere^7= U(3)\backslash U(4)$ the coinduced Poisson
bracket $(\sphere^7,\{,\})$.
The bracket on $\sphere^7$ can be written as the restriction of the following
bracket in $\Co^4$: if $z_i$, $i=1,\ldots,4$, denote complex coordinates
we let

\begin{eqnarray*}
 &\{z_i,z_j\}=z_iz_j,\quad 1\le i<j\le 4\qquad
 &\{z_i^*,z_j^*\}=-z_i^*z_j^*,\quad 1\le i<j\le 4\\
 &\{z_i,z_j^*\}=-z_iz_j^*,\quad 1\le i\ne j\le 4\qquad
 &\{z_j^*,z_j\}=\sum_{i<j}z_jz_j^* \, .
\end{eqnarray*}

More detailed information can be found in \cite{VS}.

The Lie algebra of the diagonal $SU(2)^d$ is ${\bf
su}(2)^d=\langle H_1+H_3,E_1+E_3,F_1+F_3\rangle$; using
(\ref{coboundary}) it is easy to verify that $\delta_{\Re}({\bf
su}(2)^d)\not\subset {\bf su}(2)^d\wedge \u4$ so that $SU(2)^d$ is
not a coisotropic subgroup. We then have to solve the following problem:

\medskip
{\it Does there exist any $g\in U(4)$ such that $\delta_{\Re}(\Ad_g(
{\bf su}(2)^d))\subset \Ad_g({\bf su }(2)^d)\wedge \u4$, i.e. such that
$gSU(2)^dg^{-1}$ is coisotropic ?}

\medskip
We give a positive answer to this question. By direct computation we
verify that if
$$
g=\left(\begin{array}{cccc} 1&0&0&0\cr 0&1&0&0\cr 0&0&0&1\cr
0&0&-1&0
\end{array}\right)  \in U(4) ~~~,
$$
then $SU(2)^d_g=g SU(2)^d g^{-1}$ is a coisotropic subgroup of $U(4)$.
This is not the only solution but the general problem will be studied
elsewhere.
The projection onto this subgroup is then defined by
\begin{equation}
\label{twisted_projection}
r_g(t) = \left(
\begin{array}{cccc}
\alpha & \beta & 0 & 0 \cr -\beta^* & \alpha^* & 0 & 0\cr 0 & 0 &
\alpha^* & \beta^* \cr 0 & 0 & -\beta & \alpha
\end{array}
\right)\,,~~~~~~~~ |\alpha|^2+|\beta|^2=1 \;.
\end{equation}

The right action of $SU(2)^d_g$ on $\sphere^7\simeq U(3)\backslash U(4)$
is free and defines a principal bundle on  $\sphere^4 \simeq U(3)\backslash
U(4)/SU(2)^d_g$ which is isomorphic to the Hopf bundle.
Indeed it is easy to verify that $i:\sphere^7\to\sphere^7$,
$i(z_1,z_2,z_3,z_4)=(z_1,z_2,-z_4,z_3)$ is a bundle morphism.
Nevertheless since $SU(2)^d_g$ is coisotropic on $U(4)$, thanks to
Proposition \ref{prop_Lu}, a Poisson structure is coinduced on the base
and the projection $\sphere^7\to\sphere^4$ is a Poisson map.
We call this bundle a {\it Poisson principal bundle}.

The Poisson structure can be explicitly described by the restriction of
the bracket of $\sphere^7$ to the subalgebra generated by the following
coinvariant functions
$$
R=z_1z_1^*+z_2z_2^*,\qquad
a=z_1z_4^*-z_2z_3^*,\qquad b=z_1z_3+z_2z_4\ ,
$$
which satisfy $|a|^2+|b|^2=R(1-R)$ .
Easy calculations prove that:
\begin{eqnarray*}
&\{a,R\}=-2aR,\qquad
\{b,R\}=2b R,\quad &
\{a,b\}=-3ab, \qquad
\{a,b^*\}=ab^*, \quad\\
&\{a,a^*\}=-2aa^*+2R^2,&
\quad \{b,b^*\}=4bb^*-2R \ .
\end{eqnarray*}

This Poisson algebra has clearly zero rank in $R=0$.
Let $R\ne 0$ and define $\zeta_1=a/R$, $\zeta_2=b/R$.
Geometrically we're just giving cartesian coordinates on the stereographic
projection on $\Co^2$.
Poisson brackets between these new coordinates are given by:
\begin{equation}
\label{symplectic_complex}
    \begin{array}{ll}
    \{\zeta_1,\zeta_2\}=\zeta_1\zeta_2,\quad
    &\{\zeta_1,\zeta_1^*\}=2(1+|\zeta_1|^2),
    \\ \{\zeta_1,\zeta_2^*\}=\zeta_1\zeta_2^*, &
    \{\zeta_2,\zeta_2^*\}=-2(1+|\zeta_1|^2+
    |\zeta_2^2|)\;.
    \end{array}
\end{equation}
\noindent  Such brackets define a symplectic structure on the
$4$--dimensional real space $\Re^4$ (it can be proven, in fact
that the corresponding map between cotangent and tangent bundle
has fixed maximal rank).
The covariant Poisson bracket on $\sphere^4$ has thus a very simple
foliation given by a $0$--dimensional leaf and a 4--dimensional linear space.

We summarize this discussion in the following Proposition.
\medskip
\begin{proposition}
\label{su2dg}
The embedding of $SU(2)_g^d$ into $U(4)$ defines a coisotropic subgroup.
The corresponding bundle $\sphere^7\to\sphere^4\simeq\sphere^7/SU(2)_g^d$
is a Poisson bundle.
\end{proposition}


\section{The Quantum $\sph4$}
The Hopf algebra $U_q(4)$ is generated by $\{t_{ij}\}_{ij=1}^4$ ,
$D_q^{-1}$ and the following relations (see \cite{KS}):
\begin{eqnarray*}
t_{ik} t_{jk} = q \ t_{jk} t_{ik}\;,  & {} & t_{ki} t_{kj} = q \
t_{kj} t_{ki}\;,  ~~~~~~~~~i<j \cr t_{i\ell} t_{jk} &=&  t_{jk}
t_{i\ell}\;, ~~~~~~~~~~~~~i<j,\,k<\ell\cr t_{ik} t_{j\ell} -
t_{j\ell} t_{ik} & = & (q-q^{-1}) t_{jk} t_{i\ell}\;, ~~~~~~~~~i<j,\,
k<\ell\cr D_qD_q^{-1} &=& D_q^{-1} D_q = 1\;,
\end{eqnarray*}
where $D_q =\sum_{\sigma\in
P_4}(-q)^{\ell(\sigma)}t_{\sigma(1)1}\ldots t_{\sigma(4)4}$ with
$P_4$ being the group of 4--permu\-tations, is central.
The Hopf algebra structure is
\begin{eqnarray*}
\Delta(t_{ij})=\sum_k t_{ik}\otimes t_{kj}\;,~~&~~
\Delta(D_q)=D_q\otimes D_q\;, \cr
\epsilon(t_{ij})=\delta_{ij}\;,~~&\epsilon(D_q)=1\;,
\end{eqnarray*}
$$
S(t_{ij}) = (-q)^{i-j} \sum_{\sigma\in P_3(j)} (-q)^{\ell(\sigma)}
t_{\sigma(1)1}\ldots t_{\sigma(i)i+1}\ldots t_{\sigma(4)4} ~~
D_q^{-1} \;,
$$
the compact real structure forces us to choose $q\in\Re$ and is defined by
$t^*_{ij} = S(t_{ji})$, $D_q^* = D_q^{-1}\;.$

\medskip
Let $L = {\rm Span}\{t_{j4},t_{4j},t_{44}-1\}_{j=1,2,3}$; it comes
out that $\L = U_q(4)L$ is a Hopf ideal, so that
$U_q(4)/\L$ is equivalent to $U_q(3)$ as a Hopf algebra.
Let $\ell:U_q(4)\to U_q(4)/\L\simeq U_q(3)$ be the quotient
projection.
The algebra $\sphere^7_q = {}^{\rm \ell }U_q(4)$ is generated by
$z_i = t_{4i}$, $i=1\ldots4$, with the following relations \cite{VS}:
\begin{eqnarray*}
z_i z_j\ =\ q z_j z_i  ~~(i<j) ~,~~~~~~~~~ z_j^*z_i = q z_i z_j^*
~~(i\neq j)~, \cr z^*_k z_k = z_k z^*_k + (1-q^2) \sum_{j<k} z_j
z_j^* ~,~~~~~~~~~\sum_{k=1}^4 z_k z^*_k = 1\;.
\end{eqnarray*}
The $U_q(4)$-coaction on $\sphere^7_q$ reads $\Delta(z_i)=\sum_j
z_j\otimes t_{ji}$.

\medskip
Let us now quantize the coisotropic subgroup $SU(2)^d_g$ of
Proposition \ref{su2dg}.
Motivated by the projection $r_g$ in (\ref{twisted_projection}) let us
define $\R=R\ U_q(4)$, where
\begin{eqnarray*}
R&=&{\rm Span}\{t_{13},\,t_{31},\,t_{14},\,t_{41},\,t_{24},\,t_{42},
\,t_{23},\,t_{32},\,t_{11}-t_{44},\,t_{12}+t_{43},\cr
&&\qquad\ \ t_{21}+t_{34},\,t_{22}-t_{33},\,t_{11}t_{22}-q\ t_{12}t_{21}-1\}\cr
&=&\dot{R}\oplus {\rm Span}\{
t_{11}t_{22}-q\ t_{12}t_{21}-1\}\ .
\end{eqnarray*}
It is easy to verify that $\R$ is a $\tau$--invariant, right ideal,
two sided coideal. Let $r:U_q(4)\to U_q(4)/\R$ be the
projection map.
We have the following result:
\medskip
\begin{proposition}
\label{suq2}
As a $\tau$--coalgebra $U_q(4)/\R$ is isomorphic to
$SU_q(2)$.
\end{proposition}
\smallskip
\noindent
{\it Proof}. We sketch here the main lines of the proof.
Let $A_q(N)$ be the bialgebra generated by the $\{t_{ij}\}$.
We first remark that $r(D_q)=1$ so that $U_q(4)/\R\simeq A_q(4)/R A_q(4)$.
First one can show that $A_q(4)/ \dot{R}A_q(4)\simeq A_q(2)$.
Once chosen an order in the generators $t_{ij}$ of  $A_q(4)$, a linear
basis is given by the ordered monomials in $t_{ij}$
\cite{Koelink}, so that $A_q(4)={\rm Span}\{
t_{11}^{n_{11}}\,t_{44}^{n_{44}}\,t_{12}^{n_{12}}\,t_{43}^{n_{43}}\,
t_{21}^{n_{21}}\,t_{34}^{n_{34}}\,t_{22}^{n_{22}}$
$t_{33}^{n_{33}}\}$. Making a repeated use of the following
relations for $i<k,j<l$
\begin{eqnarray*}
t_{ij}^n t_{kl}&=&t_{kl}t_{ij}^n - q^{-1}(1-q^{2n})\,t_{il}t_{kj}t_{ij}^{n-1}
\ ,\cr
t_{ij}t_{kl}^m&=&t_{kl}^mt_{ij} + q\, (1-q^{-2m})\,t_{il}t_{kj}t_{kl}^{m-1}\ ,
\end{eqnarray*}
we get that $A_q(4)/\dot{R} A_q(4)={\rm Span}\{t_{11}^{n_{11}}\,t_{12}^{n_{12}}
\,t_{21}^{n_{21}}\,t_{22}^{n_{22}}\}\simeq A_q(2)$. To show that this a
$\tau$--coalgebra isomorphism is equivalent to verify that the projection $r$
restricted to the first quadrant of $A_q(4)$ is a $\tau$--bialgebra
isomorphism.
This can be directly done by using the relations.
Finally the quotient by the quantum determinant gives $SU_q(2)$. \fidi
\medskip
\begin{remark}{\rm
The projection map $r: U_q(4)\to SU_q(2)$ is not a Hopf algebra
map as can be, for instance, explicitly verified on $r(t_{11}
t_{43})\not= r(t_{11}) r( t_{43})$.}
\end{remark}
\medskip
In the following we will denote $U_q(4)/\R$ with $SU_q(2)$, but we
have to be careful that $r$ doesn't preserve the algebra structure
but only defines a right $U_q(4)$--module structure on the quotient.

By construction $\Delta_r = (\id\otimes r)\Delta: \sphere^7_q\to
\sphere^7_q\otimes SU_q(2)$ defines an $SU_q(2)$ coaction on $\sphere^7_q$.
The space of functions on the quantum $4$--sphere $\sph4$ is the
space of coinvariants with respect to this coaction, {\it i.e.}
$\sph4=\{a\in \sphere^7_q \, | \, \Delta_r(a) = a\otimes r(1) \}$.
We describe $\sph4$ in the following proposition whose proof is postponed in
the Appendix.
\medskip
\begin{proposition}
\label{s4q}
The algebra $\sph4$ is generated by $\{a,a^*,b,b^*,R\}$, where
$a = z_1 z_4^*-z_2 z_3^*,\  b = z_1z_3 + q^{-1}z_2 z_4,
\ R = z_1 z_1^* + z_2 z_2^*$. They satisfy the following relations
$$
R a = q^{-2} a R ~,~~~ R b = q^2 b R   ~,~~~ a b = q^3 b a ~,~~~
ab^*= q^{-1} b^* a, $$
$$aa^*+q^2bb^* = R(1-q^2R),$$
$$ aa^*= q^2a^*a + (1-q^2)R^2 ~,~~~~~ b^*b =
q^4bb^* + (1-q^2)R\;.
$$
\end{proposition}
\medskip
In terms of $r_{ij}=r(t_{ij})\in SU_q(2)$, with $i,j=1,2$, the
fundamental corepresentation $\tau_f:\Co^2\to SU_q(2)\otimes
\Co^2$ is written as
$$
\tau_f\left(\begin{array}{c}e_1\cr e_2\end{array}\right)=
\left(\begin{array}{cc}r_{11}&r_{12}\cr r_{21}& r_{22}
\end{array}\right)\otimes \left(\begin{array}{c}e_1\cr
e_2\end{array}\right)\;.
$$
Let ${\cal E}=\cotens{\sphere^7_q}{\tau_f}{\Co^2} $ the associated quantum
vector bundle. Let $\langle (a_1,a_2),(b_1,b_2)\rangle= a_1
b_1^*+a_2b_2^*\in \sph4$ for $(a_1,a_2), (b_1,b_2)\in{\cal E}$ be
the hermitian structure in ${\cal E}$. Let
$$
f_1= q (z_1,z_2) ~,~~~ f_2=q(z_2^*,-qz_1^*) ~,~~~ f_3=(z_4,-z_3)
~,~~~ f_4 = q(z_3^*,q^{-1}z_4^*) \;,
$$
and $G\in M_4(\sph4)$ such that $G_{ij} = \langle f_i,
f_j\rangle$.
We then have the following description of ${\cal E}$ (see the Appendix for
the proof).
\medskip
\begin{proposition}
\label{projector}
As a $\sph4$--module ${\cal E}$ is generated by
$f_i$, $i=1\ldots 4$; it is isomorphic to $(\sph4)^4 G$
where
\begin{equation}
\label{projector_G}
G = G^2 = \left(\begin{array}{cccc} q^2 R & 0
& qa& q^2b\cr 0&q^2R&  q b^*&-q^3 a^* \cr q a^* & q b & 1-R & 0
\cr q^2 b^* & - q^3 a
& 0& 1-q^4 R
\end{array}\right)\;.
\end{equation}
\end{proposition}


\section{Unitary Representations of $\sph4$}
Let $0<q<1$.
By restriction of those of $\sphere_q^7$, see for instance
\cite{Chari_Pressley}, we obtain the following two inequivalent unitary
representations of $\sph4$.
The first is one dimensional and it is obtained as the restriction of the
counit $\epsilon$ of $U_q(4)$:
\begin{equation}
\label{rep0} \epsilon(R)=\epsilon(a)=\epsilon(b)=0\;.
\end{equation}
The second one $\sigma:\sph4\to B(\ell^2(\Na)^{\otimes 2})$
is defined by
\begin{eqnarray}
\label{irrep} \sigma(R)|n_1,n_2\rangle &=& q^{2(n_1+n_2)}
|n_1,n_2\rangle, \cr
\sigma(a)\;|n_1,n_2\rangle &=&
q^{n_1+2n_2-1}(1-q^{2n_1})^{1/2}|n_1-1,n_2\rangle,\cr
\sigma(b)\;|n_1,n_2\rangle &=&
q^{n_1+n_2}(1-q^{2(n_2+1)})^{1/2}|n_1,n_2+1\rangle \;.
\end{eqnarray}
There are no other irreducible representations
with bounded operators.
In fact let $\rho:\sph4\to B({\cal H})$ be such a representation, since
$\rho(R)$ is a bounded selfadjoint operator and $Ra=q^{-2}aR$,
$Rb^*=q^{-2}b^*R$, there exists a vector $|\lambda\rangle$ such that
$\rho(R)\,|\lambda\rangle=\lambda\, |\lambda\rangle$ and $\rho(a)|
\lambda\rangle=\rho(b^*)|\lambda\rangle=0$.
By using the relation
$a^*a+bb^*=q^{-2}R(1-R)$ we conclude that $\lambda=0$ or $\lambda=1$.
Being $\rho$ irreducible it can be verified that for $\lambda=0$ we have that
$\rho=\epsilon$ and for $\lambda=1$ we have $\rho=\sigma$.

Let us remark that such irreducible (unitary) representations are in one
to one correspondence with the leaves of the symplectic foliation of the
underlying Poisson 4--sphere: the 0--dimensional leaf corresponds to the
counit and the symplectic $\Re^4$ to the infinite dimensional representation.
The representation $\sigma$ has the following important property.
\medskip
\begin{proposition}
\label{trace_class} The operator $\sigma(x)\in
B(\ell^2(\Na)^{\otimes 2})$ is a trace class operator for each
$x\in\bar{\sph4}=\sph4/{\Co\, 1}$.
\end{proposition}
\smallskip
\noindent
{\it Proof}. Since the family of trace class operators ${\cal
I}_1$ is a $*$--ideal in the algebra of bounded operators, it is
enough to verify the proposition on the generators $\sigma(R)$,
$\sigma(a)$ and $\sigma(b)$. Indeed we have that
$\tr(\sigma(|R|))=\tr(\sigma(R))=(1-q^2)^{-2}$ and
$\tr(\sigma(|a|))=\sum_{n_1,n_2\geq
0}q^{n_1+2n_2-1}(1-q^{2n_1})^{1/2}= q^{-1}(1-q^2)^{-1}\sum_{n\geq
0}q^n(1-q^{2n})^{1/2}\leq q^{-1}(1-q^2)^{-1}\sum_{n\geq 0}q^n$
$(1-q^{2n})= (1-q)^{-1}(1-q^3)^{-1} $. Analogously $\tr(|b|)\leq
(1+q^2)(1-q)^{-1}(1-q^3)^{-1}$. \fidi
\medskip
\begin{remark}{\rm
The universal $C^*$--algebra $C(\sph4)$, defined by $\sph4$, is the norm
closure of $\sigma(\sph4)$.
By Proposition \ref{trace_class} we have that $\sigma(\sph4)
\setminus \Co 1$ is contained in the algebra $K$ of compact operators on
$B(\ell^2(\Na)^{\otimes 2})$. Using Proposition 15.16 of \cite{DorFel}
we conclude that $C(\sph4)$ is isomorphic to the unitization of compacts.

Note that, though different at an algebraic level, it is  not
possible to distinguish from their C*--algebras our 4--sphere and
the standard Podle\`s sphere $\sphere^2_q(c,0)$ in
\cite{MasNakWat,Podles}. A possible explanation for this
peculiarity stands in the fact that the space of leaves of the
underlying symplectic foliations are homeomorphic.}
\end{remark}
\medskip
Let $H=\ell^2(\Na)^{\otimes 2}\oplus \ell^2(\Na)^{\otimes 2}$ and
$\pi=\left ( \begin{array}{rr}\sigma & 0\\0 &\epsilon\end{array}
\right ),\ \gamma=\left ( \begin{array}{rr}1 & 0\\0 &-1\end{array}
\right )$. As a consequence of Proposition \ref{trace_class} we have that
$(H,\pi)$ is a 0--summable Fredholm module whose character is
$tr_\sigma=\tr\; (\sigma -\epsilon)$.
By explicit computation we have $tr_\sigma(1)=0$ and
\begin{eqnarray}
\label{trace} \tr_{\sigma}(R^{k}) =
\frac{1}{(1-q^{2k})^2}~~(k>0)~, ~~~~~
\tr_{\sigma}(a)=\tr_{\sigma}(b) = 0 ~,~~   \cr
\tr_{\sigma}(aa^*)=\frac{1}{(1-q^4)^2} ~,~~
\tr_{\sigma}(bb^*)=\frac{1}{(1-q^2)(1-q^4)}~~~~\;.
\end{eqnarray}

In the representation $\sigma$, $R$ is an invertible operator.
This suggests the possibility of realizing a {\it quantum stereographic
transformation} on a deformation of $\Co^2$, that we will denote $\Co^2_q$.
Let us define $\zeta_1=R^{-1}a$, $\zeta_2=b R^{-1}$; by direct computation
we find that they satisfy the following relations:
\begin{equation}
    \begin{array}{ll}
    \zeta_1\zeta_2= q^{-1}\zeta_2\zeta_1,\quad
    &\zeta_1\zeta_1^*=q^{-2}\zeta_1^*\zeta_1 + (1-q^2),
    \\ \zeta_1\zeta_2^*=q^{-1}\zeta_2^*\zeta_1, &
    \zeta_2\zeta_2^*=q^2\zeta_2^*\zeta_2
    -(1-q^2)(q^2+\zeta_1^*\zeta_1)~~.
    \end{array}
\end{equation}
One can verify that the algebra $\Co^2_q$ quantizes the symplectic
structure on $\Co^2$ seen in (\ref{symplectic_complex}).


\section{Chern--Connes paring of G}

Let us compute the Chern classes in cyclic homology associated
to the projector $G$.
We briefly recall some basic definitions and results from cyclic
homology, see \cite{ConnesBook} and \cite{Loday} for any details.

Let $A$ be an associative $\Co$--algebra.
Let $d_i(a_0\otimes a_1\ldots \otimes a_n )=a_0\otimes..\,a_ia_{i+1}..
\otimes a_n$, for $i=0,\ldots n-1$ and $d_n(a_0\otimes a_1\ldots \otimes
a_n)=a_na_0\otimes a_1\ldots \otimes a_{n-1}$; the Hochschild
boundary is defined as $\beta=\sum_{i=0}^n(-)^id_i$ and the
Hochschild complex is ($C_*(A),\beta$), with $C_n(A)=A^{\otimes
n+1}$.
As usual we denote Hochschild homology with $HH_*(A)$.

Let $t(a_0\otimes a_1\otimes \ldots a_n)=(-)^n a_1\otimes\ldots
a_n\otimes a_0$ be the cyclic operator and
$C^\lambda_n(A)=A^{\otimes n+1}/(1-t)A^{\otimes n+1}$.
The Connes complex is then ($C^\lambda_*(A),\beta$); its homology is denoted
as $H^\lambda_*(A)$.
For each projector $G\in M_k(A)$, {\it i.e.} $G^2=G$, the Chern class is
defined as $ch_n^\lambda(G)=\Tr[ (-)^nG^{\otimes 2n+1}]\in H^\lambda_{2n}(A)$,
where $\Tr:M_k(A)^{\otimes n}\to A^{\otimes n}$ is the
generalized trace, {\it i.e.} $\Tr[M_1\otimes \ldots \otimes M_n]=
\sum_{j}[M_1]_{j_1j_2}\otimes [M_2]_{j_2j_3}\ldots \otimes
[M_n]_{j_nj_1}$.

Let $I: A^{\otimes n+1}\to C_n^\lambda(A)$ be the projection map and let
$x\in A^{\otimes n}$ be such that $I(x)$ is a
cycle which induces $[I(x)]\in H_n^\lambda$.
Let us define the periodicity map $S: H^\lambda_n\to H^\lambda_{n-2}$, as
$$
S([I(x)]) = -\frac{1}{n(n-1)} [I(\sum_{0\leq i<j\leq n}
(-)^{i+j}d_id_j(x))]\;.
$$
There is then a long exact sequence in homology:
\begin{equation}
\label{connes_sequence}
\ldots\to HH_n(A)
\stackrel{I}{\to} H_n^\lambda(A)\stackrel{S}{\to}
H_{n-2}^\lambda \stackrel{B}\to HH_{n-1}(A)\ldots\;,
\end{equation}
where $B$ is an operator we don't need to define. With our
normalization of the Chern character, we have that for each
projector $G$
\begin{equation}
\label{S_operator} S(ch^\lambda_n(G)) =-
\frac{1}{2(2n-1)}ch^\lambda_{n-1}(G) \;.
\end{equation}

Let $G\in M_4(\sph4)$ be the projector defined in Proposition
\ref{projector}. Then
$$
ch_0^\lambda (G)= [\Tr(G)] = [2 - (1-q^2)^2\;R] \in
H_0^\lambda\;.
$$
The character $\tr_\sigma$ of the Fredholm module $(H,\pi)$ given in
(\ref{trace}) is a well defined cyclic $0$--cocycle on $\sph4$.
We then have that $\tr_\sigma(ch_0^\lambda(G))=-1$ and conclude that
$ch^\lambda_0(G)$ defines a non trivial cyclic cycle in
$H^\lambda_0(\sph4)$; using the $S$--operator (\ref{S_operator})
and Connes sequence (\ref{connes_sequence}) we conclude that
$ch_1^\lambda$ and $ch_2^\lambda$ define non trivial classes in
cyclic homology and are not Hochschild cycles.
Since $\tr_\sigma$ is the character of a Fredholm module, the integrality of
the pairing is a manifestation of the so called noncommutative index
theorem \cite{ConnesBook}.
We summarize this discussion in the following proposition.
\medskip
\begin{proposition}
\label{chern_class}
The projector $G$ defined in {\rm(\ref{projector_G})}
defines non trivial cyclic homology classes
$ch^\lambda_n(G)\in H^\lambda_{2n}$. The Chern--Connes
pairing with $\tr_\sigma$ defined in {\rm(\ref{trace})} is:
$$
\langle \tr_\sigma, G \rangle = -1  \;.
$$
\end{proposition}
\bigskip
\begin{appendix}
\noindent
{\large\bf Appendix.} {\bf Proof of Proposition \ref{s4q} and \ref{projector}}

\medskip

\noindent
To prove Proposition \ref{s4q} and \ref{projector}  we use the strategy
adopted by Nagy in \cite{Nagy}.
His argument is based on the general fact that the corepresentation theory
of compact quantum groups is ``equivalent'' to the classical one.
This equivalence is realized by a bijective map between quantum and
classical finite dimensional corepresentations: this map preserves direct
sum, tensor product and dimension.
The fact that we deal with coisotropic subgroups requires some additional care.
The projection $r$ induces a mapping $r[\rho]=(\id\otimes r)\rho$ from the
corepresentations of $U_q(4)$ into the corepresentations of $SU_q(2)$, since
$r$ is not an algebra morphism it is not obvious that this mapping preserve
the tensor product.
However in our case the following lemma can be proved:
\medskip
\begin{lemma}
\label{tens}
Let $t_f$ and $t_{f^c}$ be the fundamental and its contragredient
corepresentation of $U_q(4)$, then
$r[t_f^{\otimes r}\otimes t_{f^c}^{\otimes s}]$ is equivalent to
$r[t_f]^{\otimes r}\otimes r[t_{f^c}]^{\otimes s}$.
\end{lemma}
\smallskip
\noindent
{\it Proof.} Let $\tau_f$ be the fundamental corepresentation of $SU_q(2)$ and
$\tau_{f^c}$ its contragredient.
Let us notice that for $i,j=1,2$ we have that $r[t_{i,j}]=(\tau_f)_{ij}$,
$\ r[t_{i+2,j+2}]=q^{i-j}(\tau_{f^c})_{ij}$, $\ r[t^*_{i,j}]=(\tau_{f^c})_{ij}$
and $r[t^*_{i+2,j+2}]=q^{j-i}(\tau_f)_{ij}$.
By making use of the equivalence between $\tau_f$ and $\tau_{f^c}$ it is
easy to conclude that $r[t_f^{\otimes r}\otimes t_{f^c}^{\otimes s}]$
for $r+s=2$ is equivalent to $4\; \tau_f^{\otimes 2}$ and then to
$r[t_f]^{\otimes r}\otimes r[t_{f^c}]^{\otimes s}$.
The result for generic $r$ and $s$ is obtained by recurrence and by making
use of the right $U_q(4)$--module structure of the projection $r$. \fidi

\medskip
As a consequence the decomposition of $r[t_f^{\otimes n}]$ into irreducible
corepresentations of $SU_q(2)$ is the same as the classical one.

Let $t_n =\bigoplus_{r+s\leq n} t_{r,s}$ where
$t_{r,s}=t_f^{\otimes r}\otimes t_{f^c}^{\otimes s}$.
We denote with $C(t_n)\subset U_q(4)$ the subcoalgebra of the matrix elements
of $t_n$. We then have $U_q(4)=\bigcup_{n\in\Na}
C(t_n)$ and we define $\sphere^7_{q,n}=C(t_n)\cap \sphere_q^7$.
Obviously $\sphere^7_{q,n}$ is a $U_q(4)$--comodule with coaction
$\Delta_n=\Delta|_{\sphere^7_{q,n}}$.
From the decomposition into irreducible corepresentations
$\Delta_n=\sum_{\lambda\in I(U_q(4))} m_\lambda\; \lambda$, where
$m_\lambda\in\Na$, we get $\sphere^7_{q,n}=\doubleoplus{\lambda\in
I(U_q(4))}{j=1,\cdots, m_\lambda}
\sphere^{7\,{\lambda,\,j}}_{q,n}$.

Let $\rho: V\to SU_q(2)\otimes V$ be an irreducible $SU_q(2)$ corepresentation.
We prove the following Lemma.
\medskip
\begin{lemma}
\label{dim}
The dimension of $\cotens{\sphere^7_{q,n}}{\rho}V$ doesn't depend on $q$.
\end{lemma}
\smallskip
\noindent
{\it Proof.} Let $P_\rho:\sphere^7_q\otimes V \to \cotens{\sphere^7_q}
{\rho} V$ be the projection defined by $P_\rho(f\otimes v)=
\sum_{(f, v)}f_{(0)}h(f_{(1)}S(v_{(-1)}))\otimes v_{(0)}$, where $h$ is the
Haar measure on $SU_q(2)$. We obviously have
that $P_\rho(\sphere^7_{q,n}\otimes V)=\doubleoplus{\lambda
\in I(U_q(4))}{j=1,\cdots, m_\lambda}
P_\rho(\sphere^{7\,\lambda,\,j}_{q,n}\otimes V)$.
Then $dim\, P_\rho(\sphere^7_{q,n}\otimes V)=\sum_{\lambda\in I(U_q(4))}
m_\lambda\, m_\rho(\lambda)$, where
$m_\rho(\lambda)=dim\, P_\rho(\sphere^{7\, \lambda,\,j}_{q,n}\otimes V)$
equals the multiplicity of $\rho$ in the decomposition of
$r[\lambda]=(\id\otimes r)\,\lambda$.
Since the correspondence between classical and quantum corepresentation
preserves dimensions, the result follows. \fidi
\medskip
\noindent
{\it Proof of Proposition \ref{s4q}.} To show that $\{a,b,R\}$ are
coinvariants is a direct computation. Let $B_q\subset \sph4$ be the *--algebra
generated by those elements.
By the use of the diamond lemma the monomials
$\{a^{*\,i_1}a^{i_2}R^jb^{*\,k_1}b^{k_2}\big|\,$ $k_1k_2=0\}$ are linearly
independent and they form a basis of $B_q$.
Note that the same monomials form a basis for the polynomial functions on
the classical 4--sphere, and define a vector space isomorphism which maps
$B_{q,n}=C(t_n)\cap B_q \to P_\K(\sphere^7_{1,n})$, where
$P_\K=P_\rho$ with $\rho$ being the identity corepresentation. Using the
Lemma \ref{dim} we then have $B_q=\sph4$. \fidi
\medskip
\noindent
{\it Proof of Proposition \ref{projector}.} By a direct check it is easy to
see that $f_i$ are in ${\cal E}=\cotens{\sphere^7_q}{\tau_f}{\Co^2}$ and that
the mapping $f_i\to e_i\, G$, with $(e_i)_j=\delta_{ij}$, is a
$\sph4$--module morphism.
Since in the classical case it is clearly bijective the result follows by
repeating the same arguments of the proof of Proposition \ref{s4q} and
applying the Lemma \ref{dim} with $\rho=\tau_f$. \fidi

\end{appendix}

\bigskip
\noindent
{\bf\large Acknowledgments}

\medskip
\noindent
The authors want to thank L. Dabrowski and G. Landi for having stimulated
this work and for the useful discussions on the subject and the referee for
the constructive suggestions.
One of us (N.C.) would like to thank A.J.L. Sheu for his comments on
the paper.

\end{document}